# TACHIBANA, KILLING AND PLANARITY NUMBERS OF COMPACT RIEMANNIAN MANIFOLDS

Sergey E. Stepanov

*Department of Mathematics, Finance University under the Government of Russian Federation, Moscow 125468, Russian Federation*

Josef Mikeš

*Department of Algebra and Geometry, Palacky University, 77146 Olomouc, Czech Republic*

**Abstract.** We present definitions and properties of conformal Killing, Killing and planarity forms on a Riemannian manifold and determine Tachibana, Killing and planarity numbers as an analog of the well known Betti numbers. We state some set of conditions to characterize these numbers. Moreover, we formulate the main results on the relationship between the Betti, Tachibana, Killing and planarity numbers.

## 1. Introduction

The present paper based on our report at XVII Geometrical Seminar in Zlatibor (Republic of Serbia, September 3-8, 2012) and on our plenary lecture at AGMP-8 in Brno (Czech Republic, September 12-14, 2012). The whole paper is organized as follows. In first section, we consider a brief review of the Riemannian geometry of conformal Killing, Killing and closed conformal Killing forms. In the next paragraph of the paper we determine the Tachibana number $t_r(M)$, Killing number $k_r(M)$ and planarity number $p_r(M)$ as an analog of the well known Betti number $b_r(M)$ of $(M, g)$ for $1 \leq r \leq n-1$ and show the properties of these numbers. Finally, we state some set of conditions to characterize these numbers and formulate and prove results on the relationship between the Betti, Tachibana, Killing and planarity numbers.

## 2. A brief review of the Riemannian geometry of conformal Killing, Killing and closed conformal Killing forms

### 2.1. Conformal Killing forms

Conformal Killing forms (which are also called in the physics literature as *conformal Killing-Yano tensors*) were introduced by Tachibana (see [35]) and Kashiwada (see [13]) as a natural generalization of *conformal Killing vector fields*, which are also called infinitesimal conformal transformations (see [21]; [43], pp. 46, 47). Since then these forms found wide applications in physics related to hidden symmetries, conserved quantities, symmetry operators, or separation of variables (see [5]; [7]; [8], pp. 414, 426; [9]; [32] and etc.). At the same time, the conformal Killing forms have been studied by many geometers (see, for example, [10]; [14]; [22]; [24] and etc.).

Let $M$ be an $n$-dimensional $C^\infty$-manifold with a pseudo-Riemannian metric tensor $g$ and the Levi-Civita connection $\nabla$. A pair $(M, g)$ is called a *pseudo-Riemannian manifold*. On our fixed manifold $M$, we denote by $\Lambda^r M$ $(r = 1,...,n)$ the exterior power $\Lambda^r(T^*M)$ of cotangent bundle $T^*M$ of $M$. Hence $C^\infty \Lambda^r M$, the space of all $C^\infty$-sections of the bundle $\Lambda^r M$ over the manifold $M$, is the vector space $\Omega^r(M)$ of all exterior differential forms of degree $r$ on $M$ for $r = 1,...,n$.

An exterior differential form $\omega \in C^\infty \Lambda^r M$ is said to be *conformal Killing r-form*, or according to another terminology, to be *conformal Killing-Yano tensor field of order r* if there is an exterior differential $(r-1)$-form $\theta$ called the associated form of $\omega$ such that

$$\nabla \omega = \frac{1}{r+1} d\omega + g \wedge \theta \qquad (2.1)$$

where $d: \Omega^r(M) \to \Omega^{r+1}(M)$ is the well known exterior differential operator and

$$(g \wedge \theta)(X_0, X_1, ..., X_r) = \sum_{a=2}^{r} (-1)^a g(X_0, X_a) \theta(X_1, ..., X_{a-1}, X_{a+1}, ..., X_r)$$

for any vector fields $X_0, X_1, ..., X_r \in C^\infty(TM)$ and for all = 1, 2, ..., $n-1$. From (2.1) we obtain $\theta = \frac{1}{n-r+1} d^*\omega$ for the co-differential operator $d^*: \Omega^r(M) \to \Omega^{r-1}(M)$ formally adjoint to the exterior differential operator $d$ (see also [13]).

An important property of the system of linear differential equations (2.1) for $\omega \in \Omega^r(M)$ is that it is conformally invariant (see [5]). This implies that if $\omega$ is a conformal Killing form on a manifold with metric tensor $g$, then $\bar{\omega}$ is a conformal Killing form with the confor-mally scaled metric $\bar{g}$ where the conformal weights are related by $\bar{g} = e^{2f} g$ and $\bar{\omega} = e^{(r+1)f} \omega$.

Moreover, in [14] it was noted that the system of linear differential equations (2.1) is invariant under the Hodge duality, i.e. the skew symmetric part $(r+1)^{-1} d\omega$ of $\nabla \omega$ transforms into the divergence part $(n-p+1)^{-1} g \wedge d^* $ of $\nabla \omega$ and vice versa. This implies that the dual $*\omega$ is a conformal Killing $(n-r)$-form whenever $\omega$ is a conformal Killing $r$-form for the well known *Hodge star operator* $*: \Omega^r(M) \to \Omega^{n-r}(M)$ (see [1], p. 33; [19], p. 203-204). Therefore, if we denote by $\mathbf{T}^r(M, \mathbf{R})$ a vector space of conformal Killing $r$-forms then the following isomorphism $*: \mathbf{T}^r(M, \mathbf{R}) \to \mathbf{T}^{n-r}(M, \mathbf{R})$ holds (see the prove of this proposition in [14]).

A system of linear differential equations (2.1) was analyzed in [13]. In particular, the author of this paper has studied the integrability conditions of these equations on Riemannian manifolds (M, g). They have proved that if the system of the equations (2.1) is a totally-integrable system of differential equations in a neighborhood of an arbitrary point of (M, g) then (M, g) is a (locally) conformal flat Riemannian manifold. On the other hand, the following equality

$$\dim \mathbf{T}^r(M, \mathbf{R}) \leq \frac{(n+1)!}{r!(n-r+1)!}$$

holds for all $r = 1, \ldots, n-1$. The equality attained on a closed conformally flat Riemannian manifold (M, g) (see [22], [33]).

Moreover, it is well known that the Lie algebra of the group $C(M, g)$ of conformal transformations of an $n$-dimensional connected Riemannian manifold has a "local dimension" not greater than ½ $(n+1)(n+2)$ and is implemented as a space of conformal Killing vector fields (see [42]). This property follows from the equality in the case $r = 1$.

*2.2. Killing forms*

Two special subclasses of conformal Killing forms are of particular interest: *Killing forms*

are those with zero divergence part of (2.1) and *closed conformal Killing forms* with vanishing skew symmetric part of (2.1).

The concept of a Killing forms (which are also called in the physics literature as *Killing-Yano tensors*) was introduced as a generalization of *Killing vector fields* (see [19], p. 23; [21]; [43], pp. 43, 44). The Killing forms have been studied by many authors and, in particular, the equations of Killing forms has been studied intensively in the physics literature in connection with its role generating quadratic first integrals of the geodesic equation (see, for example, [4]; [6]; [11]; [15]; [26]; [37]; [43] and etc.).

As it is well known, an exterior differential form $\omega \in \Omega^r(M)$ is said to be *Killing differential r-form*, or according to another terminology, to be *Killing-Yano tensor field of order r* if it satisfies (see [43], p. 68)

$$\nabla \omega = \frac{1}{r+1} d\omega. \qquad (2.2)$$

for any $r = 1, 2, \ldots, n - 1$. It is clear, that a Killing $r$-form is a coclosed conformal Killing $r$-form. Therefore the vector space $\mathbf{K}^r(M, \mathbf{R})$ of co-closed conformal Killing $r$-forms is defied as $\mathbf{K}^r(M, \mathbf{R}) = \mathbf{T}^r(M, \mathbf{R}) \cap \mathbf{F}^r(M, \mathbf{R})$ where $\mathbf{F}^r(M, \mathbf{R})$ is a vector space of coclosed $r$-forms.

An important property of the system of linear differential equations (2.2) is that it is projectively invariant (see [26]; [28]). This implies that if $\omega$ is a Killing form on a manifold with metric tensor $g$, then $\omega' = e^{(r+1)f} \omega$ for $f = \frac{1}{2(n+1)} \ln \left| \frac{\det g'}{\det g} \right|$ is a Killing form with the projectively equivalent metric tensor $g'$ (see [17]).

A system of linear differential equations (2.2) was analyzed in [37]. In particular, the authors of this paper have studied the integrability conditions of these equations on Riemannian manifolds $(M, g)$. They have proved that if the system of the equations (2.2) is a totally-integrable system of differential equations in a neighborhood of an arbitrary point of $(M, g)$ then $(M, g)$ is a (locally) manifold of constant curvature.

On the other hand, we have proved in [26] that The components of an arbitrary Killing $r$-form $(1 \leq r \leq n - 1)$ on a pseudo-Riemannian manifold $(M, g)$ of constant curvature have the following local structure

$$\omega_{i_1 i_2 \ldots i_p} = e^{(p+1)f}\left(C_{i_0 i_1 i_2 \ldots i_p} x^{i_0} + C_{i_1 i_2 \ldots i_p}\right) \quad (2.3)$$

in a neighborhood $U$ with local coordinates $x^1, \ldots, x^n$ where $C_{i_0 i_1 \ldots i_r}$ и $C_{i_1 \ldots i_r}$ are any skew symmetric constants and $f = (2n+2)^{-1} \ln|\det g|$. Therefore we can formulate the proposition as a corollary of these theorems where we can determine a sharp upper bound on the dimension of the space of Killing $r$-forms in the form

$$\dim \mathbf{K}^r(M, \mathbf{R}) \leq \frac{(n+1)!}{(r+1)!(n-r)!}$$

for all $r = 1, \ldots, n-1$. The equality attained on closed Riemannian manifold $(M, g)$ with constant curvature.

Moreover, it is well known that the Lie algebra of the group $I(M, g)$ of motions has a "local dimension" not greater than ½ $n(n+1)$ and realized as a vector space of Killing vector fields (see, for example, [42]). If we suppose that $r = 1$ then this property follows from the above proposition.

## 2.3. Closed conformal Killing forms

The concept of a closed conformal Killing exterior differential $r$-form was introduced as a generalization of a *closed conformal vector field* (see [20]). The closed conformal Killing forms (which are also called as closed conformal Killing-Yano tensors) have been studied by many authors (see, for example, [8], p.426; [12]; [25] and etc.). In particular, closed conformal Killing forms (which are also called in the physics literature as *principal conformal Killing-Yano tensors*) have an important role for the study of hidden symmetries of higher-dimensional black holes (see about it in [8] and [16]) and electrodynamics in the general relativity theory (see [23]; [32]).

An exterior differential form $\omega \in \Omega^r(M)$ is said to be *closed conformal Killing differential r-form*, or according to another terminology, to be *closed conformal Killing-Yano tensor field of order r* if there is an exterior differential $(r-1)$-form $\theta$ called the *associated form* of $\omega$ such that

$$\nabla \omega = g \wedge \theta \quad (2.3)$$

for any $r = 1, 2, \ldots, n - 1$. Therefore the space $\mathbf{F}^r(M, \mathbf{R})$ of closed conformal Killing $r$-forms is define as $\mathbf{F}^r(M, \mathbf{R}) = \mathbf{D}^r(M, \mathbf{R}) \cap \mathbf{T}^r(M, \mathbf{R})$ where $\mathbf{D}^r(M, \mathbf{R})$ is a vector space of closed $r$-forms.

An important property of two systems of linear differential equations (2.2) and (2.3) is that these systems transform into each other under the Hodge duality (see [25]). This implies that the dual $*\omega$ is a conformal Killing $(n - r)$-form whenever $\omega$ is a Killing $r$-form. The converse proposition also holds.

Two vector spaces $\mathbf{P}^r(M, \mathbf{R})$ and $\mathbf{K}^{n-r}(M, \mathbf{R})$ transform into each other under the Hodge duality. Namely, we have proved in [25] that the isomorphism $*\colon \mathbf{P}^r(M, \mathbf{R}) \to \mathbf{K}^{n-r}(M, \mathbf{R})$ holds. Therefore, we have

$$\dim \mathbf{P}^r(M, \mathbf{R}) \leq \frac{(n+1)!}{r!(n-r+1)!}$$

for all $r = 1, \ldots, n - 1$.

**Remark.** Let $(M, g)$ be an $n$-dimensional connected Riemannian manifold with constant curvature. There is the sum $\omega' + \omega''$ of an arbitrary Killing $r$-form $\omega'$ and a closed conformal Killing $r$-form $\omega''$ is a conformal Killing $r$-form for all $r = 1, \ldots, n - 1$ (see [13]). The converse is true. Using this fact, we have determined a sharp upper bound on the dimension of the space of conformal Killing $r$-forms.

We denote by $\mathbf{C}^r(M, \mathbf{R})$ the subspace of $\Omega^r(M)$ which consists of parallel $r$-forms, which are also called covariant constant $r$-forms. It is clear, that

$$\mathbf{C}^r(M, \mathbf{R}) = \mathbf{K}^r(M, \mathbf{R}) \cap \mathbf{P}^r(M, \mathbf{R}) \subset \mathbf{H}^r(M, \mathbf{R})$$

where $\mathbf{H}^r(M, \mathbf{R}) = \mathbf{F}^r(M, \mathbf{R}) \cap \mathbf{D}^r(M, \mathbf{R})$ is a vector space of harmonic $r$-forms (see [19], pp. 205-212; [43], pp. 10, 64).

Finally, using definitions and propositions, which we have formulated above, we can obtain a 3D-diagram of the rough classification of differential $r$-forms defined on a Riemannian manifold $(M, g)$ where, for instance, the arrow $\mathbf{F}^r(M, \mathbf{R}) \to \mathbf{K}^r(M, \mathbf{R})$ means that the vector space $\mathbf{K}^r(M, \mathbf{R})$ is a subspace of $\mathbf{F}^r(M, \mathbf{R})$.

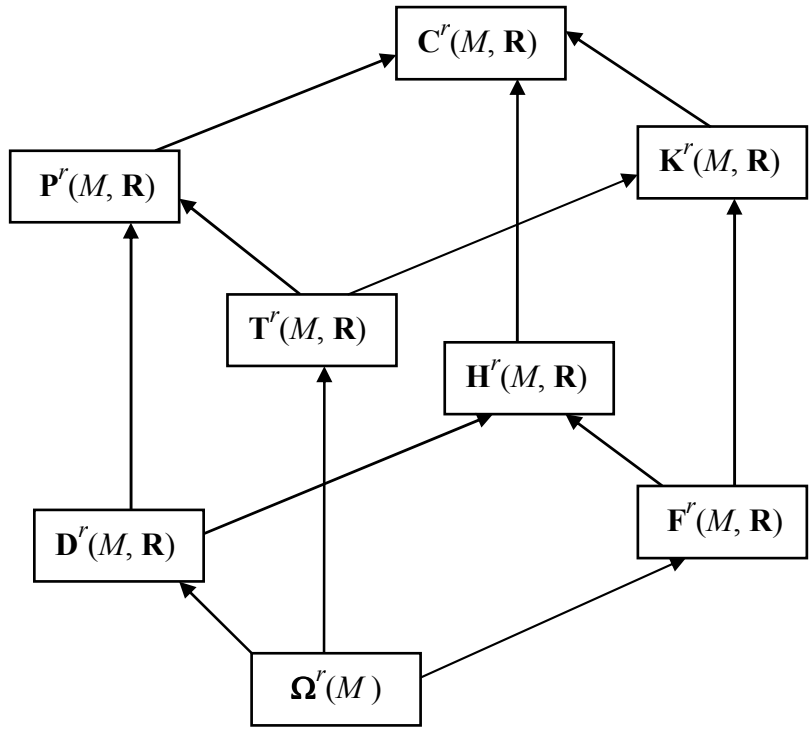

## 3. Tachibana, Killing and planarity numbers

*3.1. Rough Laplacians on differential forms*

Throughout we let $(M, g)$ be an $n$-dimensional compact and oriented Riemannian manifold with the global scalar product (see [19], p. 203)

$$\langle \omega, \omega' \rangle = \int_M g(\omega, \omega') d\text{vol} \qquad (3.1)$$

for an arbitrary $\omega, \omega' \in \Omega^r(M)$. In this case for any natural linear first-order differential operator $D$ the inner product structures on $\Omega^r(M)$ allow us to define the formal adjoint operator $D^*$ to $D$ (see [1], p. 460). Then the first order natural Riemannian (with respect to isometric diffeomorphisms) operator $D$ we can define the *rough Laplacian* $D^* \circ D$ (see [3], pp. 316-317). The first simple property of these operators comes from the fact that any rough Laplacian is a nonnegative elliptic second order linear differential operator (see [3], pp. 314-315) with a finite dimensional kernel (see also [1], pp. 461-463). For other properties of these operators, see [2]; [3]; [19]; [40]; [43] and etc.

*3.2. Definitions and properties of Tachibana, Killing and planarity numbers*

More than thirty years ago, Bourguignon (see [3]) considered the space of natural Riemannian (with respect to isometric diffeomorphisms) first-order differential operators on $\Omega^r(M)$ with values in the space of homogeneous tensors on $M$. He proved the existence of a basis of this space which consists of three operators $\{D_1, D_2, D_3\}$ where $D_1 := d$ and $D_2 := d^*$. As for the third operator $D_3$, Bourguignon said that $D_3$ does not have any geometric interpretation for $r > 1$. It was also pointed out that in the case $r = 1$ the kernel of $D_3$ consists of infinitesimal conformal transformations of $M$.

By way of specification of Bourguignon's result, we showed that (see [23]; [34])

$$D_1 = \frac{1}{p+1} d; \quad D_2 = \frac{1}{n-p+1} g \wedge d^*; \quad D_3 = \nabla - \frac{1}{p+1} d - \frac{1}{n-p+1} g \wedge d^*$$

and proved that the kernel of the third basis operator $D_3$ consists of conformal Killing $r$-forms. After that we showed in [25] that the formal adjoint operator $D_3^*$ to the natural operator $D_3$ determined by the formula $D_3^* = \nabla^* - \frac{1}{r+1} d^* - \frac{1}{n-r+1} d \circ trace$ and defined the second order differential operator

$$\Box = D_3^* \circ D_3 = \frac{1}{r(r+1)} \left( \overline{\Delta} - \frac{1}{r+1} d^* \circ d - \frac{1}{n-r+1} d \circ d^* \right)$$

which is a rough Laplacian. Here, the symbol $\overline{\Delta}$ denotes the *Bochner rough Laplacian* $\overline{\Delta} := \nabla^* \circ \nabla$ where $\nabla^*$ is the formal adjoint operator to $\nabla$ (see [1], p. 52; [2], p. 377). In particular, for $n = 2r$ we have

$$D_3^* \circ D_3 = \frac{1}{r(r+1)^2} \left( (r+1)\overline{\Delta} - \Delta \right)$$

where $\Delta := dd^* + d^*d$ is the *Laplacian on forms*, also called the *Hodge Laplacian* (see [1], p. 55; [19], p. 204; [43], p. 10).

We proved the following three propositions (see [25]):

$$\omega \in \mathbf{T}^r(M, \mathbb{R}) \iff \omega \in Ker(D_3^* D_3);$$

$$\omega \in \mathbf{K}^r(M, \mathbb{R}) \iff \omega \in Ker\, D_3^* D_3 \cap Ker\, D_2;$$

$$\omega \in \mathbf{P}^r(M,\mathbf{R}) \iff \omega \in \operatorname{Ker} D_3^* D_3 \cap \operatorname{Ker} D_1.$$

It is known that the kernel of the rough Laplacian $\square = D_3^* \circ D_3$ has a finite dimension and therefore we concluded that (see [24])

$$\dim \mathbf{T}^r(M,\mathbf{R}) = \dim_{\mathbf{R}}\left(\operatorname{Ker} D_3^* D_3\right) = t_r(M) < \infty;$$

$$\dim \mathbf{K}^r(M,\mathbf{R}) = k_r(M) < \infty; \quad \dim \mathbf{P}^r(M,\mathbf{R}) = p_r(M) < \infty$$

on a compact Riemannian manifold $(M, g)$. The numbers $t_r(M)$, $k_r(M)$ and $p_r(M)$ we have called the *Tachibana number*, the *Killing number* and the *planarity number* of a compact Riemannian manifold $(M, g)$, respectively (see [26]).

**Remark**. From the Hodge theory (see [19], pp. 202-212), we know that in a compact Riemannian manifold $(M, g)$, the number of linearly independent (with constant real coefficients) harmonic $r$-forms is equal to the Betti number $b_r(M)$ of the manifold $(M, g)$ i.e. $\dim \mathbf{H}^r(M,\mathbf{R}) = b_r(M)$. On the other hand (see [19], p. 205 ), if $(M, g)$ is a compact Riemannian manifold then we have the following definition $\mathbf{H}^r(M,\mathbf{R}) = \{\omega \in \Omega^r(M): \Delta\omega = 0\}$. This show that $\dim \mathbf{H}^r(M,\mathbf{R}) = \dim_{\mathbf{R}} \operatorname{Ker} \Delta = b_r(M) < \infty$. Thus, anyone can conclude that our definition the Tachibana numbers is an analogous to the definition of the Betti numbers.

The numbers $t_r(M)$, $k_r(M)$ and $p_r(M)$ satisfy the following duality properties $t_r(M) = t_{n-r}(M)$ and $p_r(M) = k_{n-r}(M)$ for all $r = 1, \ldots, n-1$. These equalities are analogue of the Poincare duality for the Betti numbers $b_r(M) = b_{n-r}(M)$ and corollaries of the following isomorphisms (see [14]; [23]; [27])

$$*: \mathbf{T}^r(M,\mathbf{R}) \to \mathbf{T}^{n-r}(M,\mathbf{R}); \quad *: \mathbf{P}^r(M,\mathbf{R}) \to \mathbf{K}^{n-r}(M,\mathbf{R}).$$

Moreover, the Tachibana numbers $t_r(M)$ are conformal scalar invariants, the Killing number $k_r(M)$ and the planarity number $p_r(M)$ are projective scalar invariants of a Riemannian manifold $(M, g)$ for all $r = 1, \ldots, n-1$. In addition, we remark that the first proposition is a corollary of conformal invariance of conformal Killing $r$-forms (see the propositions stated above and [5]). The second proposition of our theorem is a corollary of projective invariance of closed and co-closed conformal Killing $r$-forms (see the propositions stated above and [28]).

*3.3. On the existence of Tachibana, Killing and planarity numbers*

Let $(M, g)$ be an $n$-dimensional connected compact Riemannian manifold then the Tachibana number $t_r(M)$, the Killing number $k_r(M)$ and the planarity number $p_r(M)$ of $(M, g)$ satisfy the following inequalities

$$0 \leq t_r(M) \leq \frac{(n+2)!}{(r+1)!(n-r+1)!} \qquad 0 \leq k_r(M) \leq \frac{(n+1)!}{(r+1)!(n-r)!} \qquad 0 \leq p_r(M) \leq \frac{(n+1)!}{r!(n-r+1)!}$$

for all $r = 1, \ldots, n - 1$. Moreover, any one of two numbers $p_r(M)$ and $k_r(M)$ is maximal if and only if $(M, g)$ is Riemannian manifold with positive constant curvature. In addition, the Tachibana number $t_r(M)$ is maximal if and only if $(M, g)$ is conformal flat Riemannian manifold (see the theorems stated above and [13]; [17]; [22]).

If $\Re: \Omega^2(M) \to \Omega^2(M)$ is the well known *Riemannian curvature operator* (see, for example, [19], p. 36) then we can formulate the *vanishing theorem* of Tachibana numbers (see [24]). If $(M, g)$ is a compact and oriented Riemannian manifold with $\Re \leq 0$ then

$$t_r(M) \leq \frac{n!}{r!(n-r)!} = t_r(T^n)$$

where $T^n$ is a flat Riemannian $n$-torus and if $\Re < 0$ somewhere, then $t_r(M) = 0$ for $r = 1, \ldots, n - 1$. As a consequence, we have $k_r(M) = 0$ and $p_r(M) = 0$ for an arbitrary compact and oriented Riemannian manifold with $\Re < 0$. This proposition is an analogue of the "vanishing theorem" of the Betti numbers (see [19], p. 212).

On the other hand, we have investigated the existence of compact Riemannian manifolds with non-zero Tachibana numbers in [24].

## 4. Some relationships between Tachibana, Killing, planarity and Betti numbers

*4.1. Theorem of the decomposition the Tachibana numbers*

It is known that on a non-compact Riemannian manifolds $(M, g)$ of non-zero constant curvature $C$ the following decomposition holds (see [13])

$$\mathbf{T}^r(M, \mathbf{R}) = \mathbf{P}^r(M, \mathbf{R}) \oplus \mathbf{K}^r(M, \mathbf{R}) \tag{4.1}$$

In addition, if $(M, g)$ is compact then $C > 0$ and the decomposition (4.1) is orthogonal and global defined on $(M, g)$ (see [29]). In this case from (4.1) we obtain that the Tachibana number $t_r(M) = p_r(M) + k_r(M)$. Moreover, if $(M, g)$ is a $2r$-dimensional compact conformally flat Riemannian manifold with the positive constant scalar curvature $s$ then the orthogonal decomposition holds also (see [36]) and therefore $t_r(M) = p_r(M) + k_r(M)$. On the other hand, we have proved (see [24]) that if $(M, g)$ has the positive-define curvature operator $\mathfrak{R}$, then Tachibana numbers $t_1(M), \ldots, t_{n-1}(M)$ are decomposed as follows $t_r(M) = p_r(M) + k_r(M)$ for $r = 1,.., n-1$. In addition, we remark that the Betti numbers $b_1(M), \ldots, b_{n-1}(M)$ are zero in each of the three cases.

We shall formulate and prove the following theorem which is a generalization of these results.

**Theorem 1.** *Let $(M, g)$ be an n-dimensional ($n \geq 2$) compact and oriented Riemannian manifold with the Betti number $b_r(M) = 0$ and Tachibana number $t_r(M)$ satisfies the inequalities $t_r(M) \geq k_r(M) > 0$ then $t_r(M) - k_r(M) = p_r(M)$ for any $r = 1,.., n-1$.*

**Proof.** It well known (see [19], p. 205) that for a compact oriented Riemannian manifold $(M, g)$, the subspaces Im $d$ and Ker $d^*$ of $\Omega^r(M)$ are orthogonal complements with respect to the scalar product (3.1). Furthermore, the orthogonal complement of subspace Im $d$ in Ker $d$ coincides with the space Ker $d \cap$ Ker $d^*$ (see, for example, [18]; [20]). Therefore, in the case when $b_r(M) = 0$, the vector spaces Im $d$ and Ker $d$ coincide due to the absence of harmonic $r$-forms on $(M, g)$.

Now we consider the vector space $\mathbf{T}^r(M, \mathbf{R}) = \{\omega \in \Omega^r(M) | \omega \in \text{Ker } D_3^* D_3 \}$ of conformally Killing $r$-forms on a compact and oriented Riemannian manifold $(M, g)$. This vector space is finite dimensional and is endowed with the scalar product (3.1). Next, we suppose that the vector space $\mathbf{T}^r(M, \mathbf{R})$ has the nondegenerate subspace $\mathbf{K}^r(M, \mathbf{R}) = \{\omega \in \Omega^r(M) | \omega \in \text{Ker } D_3^* D_3 \cap \text{Im} d^*\}$ of coclosed conformally Killing $r$-forms. It is clear from the discussion above that the orthogonal complement to the vector subspace $\mathbf{K}^r(M, \mathbf{R})$ is the vector subspace $\mathbf{P}^r(M, \mathbf{R}) = \{\omega \in \Omega^r(M) | \omega \in \text{Ker } D_3^* D_3 \cap \text{Im} d\}$ of exact conformally Killing $r$-forms, i.e. $\mathbf{T}^r(M, \mathbf{R}) = \mathbf{P}^r(M, \mathbf{R}) \oplus \mathbf{K}^r(M, \mathbf{R})$, that concludes the equality $t_r(M) = p_r(M) + k_r(M)$ for $r = 1,.., n-1$. This completes the proof.

In particular, for $r = 1$ we have the following corollary.

**Corollary 1**. *Let $(M, g)$ be an n-dimensional $(n \geq 2)$ compact and oriented Riemannian manifold with the first Betti number $b_1(M) = 0$. If the first Tachibana number $t_1(M)$ and the first Killing number $k_1(M)$ satisfy the inequalities $t_1(M) > k_1(M) > 0$ then $(M, g)$ is conformal to a sphere $S^n$ in $(n+1)$-dimensional Euclidian space.*

**Proof.** Let $(M, g)$ be an $n$-dimensional $(n \geq 2)$ compact and oriented Riemannian manifold with the first Betti number $b_1(M) = 0$ and Tachibana number $t_1(M)$ satisfies the inequalities $t_1(M) > k_1(M) > 0$. Then from Theorem 1 we conclude that $p_1(M) \neq 0$ and hence there is a conformal Killing gradient 1-form $\nabla f$ for some smooth scalar function $f$ such that $\nabla^2 f = n^{-1} \Delta f \cdot g$. In addition, we recall that (see [38]) the complete Riemannian manifold $(M, g)$ of dimension $n \geq 2$ is conformal to a sphere $S^n$ in $(n + 1)$-dimensional Euclidian space if on $(M, g)$ exists a non-constant smooth function $f$ satisfying the differential equation $\nabla^2 f = n^{-1} \Delta f \cdot g$. Therefore a Riemannian manifold $(M, g)$ that satisfies the conditions of the Corollary 1 is conformal to a sphere $S^n$ in $(n + 1)$-dimensional Euclidian space.

**Remark**. Any conformal Killing vector field without zeros on $(M, g)$ can be made a Killing vector field by a conformal change of the metric $g$. Such a field is called *inessential*, otherwise it is *essential*. If we assume that $(M, g)$ is a compact Riemannian manifold admitting an essential conformal Killing vector field, then $(M, g)$ is conformally diffeomorphic with the Euclidian sphere $S^n$ (see about it [21]). Our corollary is an analogue of this proposition. Moreover, we make the assumption here that the requirement of the vanishing of the Betti numbers in Theorem 1 and Corollary 1 is superfluous.

*4.2. The relationships between Betti and Tachibana numbers*

Next we shall prove four propositions about the relationships between Betti and Tachibana numbers of a compact Riemannian manifold as corollaries of the vanishing theorems that we formulated above.

Firstly, we consider Betti and Tachibana numbers of a conformally flat Riemannian manifold and prove the following result.

**Theorem 2**. *Let (M, g) be an n-dimensional ($n \geq 3$) compact conformally flat Riemannian manifold with the non-degenerate Ricci tensor Ric then $t_r(M) \cdot b_h(M) = 0$ for all $r, h = 1, ..., n-1$.*

**Proof.** Let $\det(Ric) \neq 0$ at every point a compact Riemannian manifold (M, g) then one the two following conditions is satisfies $Ric > 0$ or $Ric < 0$. If we suppose that (M, g) is conformally flat Riemannian manifold of dimension $n \geq 3$, then in the first case we have $\Re > 0$ and the Betti numbers $b_1(M) = ... = b_{n-1}(M) = 0$ because there are no non-zero harmonic r-forms for all $r = 1, ..., n-1$ (see [4], pp. 78-79). On the other hand, in the second case we have $\Re < 0$ and Tachibana numbers $t_1(M) = ... = t_{n-1}(M) = 0$ because there are no non-zero conformal Killing r-forms for all $r = 1, ..., n-1$ (see [30]). Therefore we conclude that $t_r(M) \cdot b_h(M) = 0$ for all $r, h = 1, ..., n-1$, which finishes the proof of the theorem.

A Riemannian manifold of constant sectional curvature is an example of a conformally flat Riemannian manifold. In the following theorem we consider this manifold.

**Theorem 3**. *Let (M, g) be an n-dimensional compact and oriented Riemannian manifold with constant sectional curvature C.*

1) *If there exists the non-zero Tachibana number $t_r(M)$ for $1 \leq r \leq n - 1$ then $C > 0$ and Betti numbers $b_1(M), ..., b_{n-1}(M)$ are equal to zero.*

2) *If there exists the non-zero Betti number $b_r(M)$ for $1 \leq r \leq n - 1$ then $C < 0$ and Tachibana numbers $b_1(M), ..., b_{n-1}(M)$ are equal to zero.*

3) *If $C = 0$ then $b_r(M) = t_r(M) = \dfrac{n!}{r!(n-r)!}$ for $1 \leq r \leq n - 1$.*

**Proof.** On a compact and oriented Riemannian manifold (M, g) with constant sectional curvature $C \neq 0$ an arbitrary non-zero conformal Killing r-form $\omega$ satisfies the following equation (see [30])

$$\frac{r}{(n-r)(r+1)}\langle d\omega, d\omega \rangle + \frac{1}{n-r+1}\langle d^*\omega, d^*\omega \rangle = (n-r)C\langle \omega, \omega \rangle,$$

from which we conclude that $C > 0$. We recall that if $(M, g)$ is a Riemannian manifold of constant sectional curvature $C > 0$ then $\Re > 0$ and hence $b_1(M) = ... = b_{n-1}(M) = 0$.

On the another hand, an arbitrary non-zero harmonic $r$-form $\omega'$ on a compact Riemannian manifold of non-zero constant sectional curvature $C$ satisfies the following equation (see [4], p. 77)

$$-\langle \nabla \omega', \nabla \omega' \rangle = (n-r)C \langle \omega', \omega' \rangle,$$

from which we conclude that $C < 0$ and hence $\Re < 0$. In this case the following identities $t_1(M) = ... = t_{n-1}(M) = 0$ hold.

Finally, if $C = 0$ we have a compact flat Riemannian manifold $(M, g)$. In this case there are $\dfrac{n!}{r!(n-r)!}$ independent harmonic and conformal Killing $r$-forms. This proposition follows from the fact that we may choose a local coordinate system $x^1, ..., x^n$ in which

$$\nabla_{X_k} \omega_{i_1 i_2 ... i_n} = \dfrac{\partial \omega_{i_1 i_2 ... i_n}}{\partial x^k} = 0 \text{ for } \omega_{i_1 i_2 ... i_n} = \omega(X_{i_1}, X_{i_2}, ..., X_{i_r}) \text{ and } X_k = \dfrac{\partial}{\partial x^k}$$ (see [39], pp. 44-45; [4], p. 77). It is obviously, that any parallel form is a conformal and harmonic form simultaneously. Hence in the third case we have $b_r(M) = t_r(M) = \dfrac{n!}{r!(n-r)!}$. This finishes the proof.

The above theorem now implies

**Corollary 2**. *Let $(M, g)$ be a compact four dimensional Riemannian manifold with the vanishing Ricci tensor and the nonzero first Betti number $b_1(M)$ then $t_1(M) = t_3(M) = 4$ and $t_2(M) = 6$.*

**Proof**. It well known that every compact four dimensional Riemannian manifold with the vanishing Ricci tensor and the nonzero first Betti number is flat (see [41]). In this case from the above theorem we obtain $b_1(M) = b_3(M) = t_1(M) = t_3(M) = 4$ and $b_2(M) = t_2(M) = 6$.

As a generalization of this corollary, we derive the next result.

**Theorem 4**. *Let $(M, g)$ be a compact $n$-dimensional Riemannian manifold with the positive semi-definite Ricci tensor Ric and the nonzero first Betti number $b_1(M)$ then*

$$\frac{h!}{r!(h-r)!} \leq t_r(M) \leq \frac{n!}{r!(n-r)!}$$

*for $1 \leq r < h = b_1(M) \leq n$.*

**Proof**. Let $(M, g)$ be a compact $n$-dimensional Riemannian manifold and $b_1(M) = h \neq 0$ then there are $h$ linearly independent nonzero harmonic 1-forms $\omega_1,...,\omega_h$. We denote by $X_1,...,X_h$ the dual harmonic vector fields, i.e. $\omega_a(Y) = g(Y, X_a)$ for $a = 1,...,h$. If we suppose that the Ricci tensor Ric is positive semi-definite then $X_1,...,X_h$ are parallel (see [4], p. 37), i.e. $\nabla X_1 = 0,...,\nabla X_h = 0$. From these equalities we conclude that $\omega_1,...,\omega_h$ are parallel too. Then the $\theta_{i_1 i_2 ... i_r} = \omega_{i_1} \wedge \omega_{i_2} \wedge ... \wedge \omega_{i_r}$ for $1 \leq i_1 < ... < i_r \leq h$ are parallel and thus conformal Killing because an arbitrary parallel $r$-form $\omega$ is a conformal Killing $r$-form for $r = 1, ..., n-1$. Since these $\theta_{i_1 i_2 ... i_r}$ are $\frac{h!}{r!(h-r)!}$ linearly independent conformal Killing $r$-forms then we have $t_r(M) \geq \frac{h!}{r!(h-r)!}$ for $r < h$. Finally, if $b_1(M) = n$, then $(M, g)$ has a parallel frame. The curvature tensor of $(M, g)$ vanishes in this frame, so $(M, g)$ is flat and hence $t_r(M) = \frac{n!}{r!(n-r)!}$. This concludes the proof of the theorem.

**Acknowledgments.** Mikeš J. is supported by grant P201/11/0356 of Czech Science Foundation. Stepanov S. appreciates the hospitality of the University of Olomouc (Czech Republic) during the fall semester of 2012. We have benefited from discussion with participants of XVII Geometrical Seminar in Zlatibor (Republic of Serbia) and International conference AGMP-8 in Brno (Czech Republic).